\documentclass[a4paper,11pt]{amsart}
\usepackage{amssymb}

\input xypic
\xyoption{all}
\xyoption{arc}

\newtheorem{theo}{Theorem}[section]

\newtheorem{prop}[theo]{Proposition}

\newtheorem{lem}[theo]{Lemma}

\newtheorem{coro}[theo]{Corollary}

\def\remark#1{{\refstepcounter{theo}\label{#1}\noindent\sc Remark  
\arabic{section}.\arabic{theo} - }}

\def\equat{\refstepcounter{theo}$$~}
\def\endequat{\leqno{\boldsymbol{(\arabic{section}.\arabic{theo})}}~$$}

\newcounter{numero}[section]

\def\NM{{\mathbb{N}}}

\def\ZM{{\mathbb{Z}}}

\def\a{\alpha}
\def\b{\beta}

\def\g{\gamma}
\def\G{\Gamma}
\def\d{\delta}
\def\D{\Delta}
\def\e{\varepsilon}

\def\l{\lambda}

\def\r{\rho}
\def\s{\sigma}

\def\th{\theta}

\def\t{\tau}
\def\x{\xi}

\def\alpb{{\boldsymbol{\alpha}}}

\def\CC{{\mathcal{C}}}
\def\DC{{\mathcal{D}}}

\def\HC{{\mathcal{H}}}

\def\LC{{\mathcal{L}}}

\def\RC{{\mathcal{R}}}

\def\TC{{\mathcal{T}}}

\def\SG{{\mathfrak S}}
\def\ab{{\mathbf a}}
\def\HCB{{\boldsymbol{\HC}}}

\def\aba{{\bar{a}}}

\def\hov{{\overline{h}}}

\def\Cov{{\overline{C}}}

\def\alpb{{\boldsymbol{\alpha}}}

\def\deg{\mathop{\mathrm{deg}}\nolimits}
\def\Oplus{\mathop{\oplus}}
\def\val{\mathop{\mathrm{val}}\nolimits}

\def\to{\rightarrow}

\def\incl{\hspace{0.05cm}{\subset}\hspace{0.05cm}}

\def\vide{\varnothing}

\def\DS{\displaystyle}
\def\SS{\scriptstyle}

\def\fin{~$\SS \blacksquare$}
\def\finl{~$\SS \square$}
\def\el{\ell}
\def\infspe{\hspace{0.1em}\mathop{\preccurlyeq}\nolimits\hspace{0.1em}}

\def\lexp#1#2{\kern\scriptspace\vphantom{#2}^{#1}\kern-\scriptspace#2}
\def\le{\hspace{0.1em}\mathop{\leqslant}\nolimits\hspace{0.1em}}
\def\ge{\hspace{0.1em}\mathop{\geqslant}\nolimits\hspace{0.1em}}
\def\notle{\hspace{0.1em}\mathop{\not\leqslant}\nolimits\hspace{0.1em}}

\mathchardef\lllllll="3278
\def\SEC{$\lllllll$}
\mathchardef\inferieur="321E
\mathchardef\superieur="321F

\def\eqna{\begin{eqnarray*}}
\def\endeqna{\end{eqnarray*}}

\def\proof{\noindent{\sc{Proof}~-} }

\def\itemth#1{\item[${\mathrm{(#1)}}$]}

\def\lel{\hspace{0.1em}\mathop{\leqslant}_\LC\nolimits\hspace{0.1em}}
\def\notlel{\hspace{0.1em}\mathop{\not\leqslant}_\LC\nolimits\hspace{0.1em}}

\def\ler{\hspace{0.1em}\mathop{\leqslant}_\RC\nolimits\hspace{0.1em}}

\def\lelr{\hspace{0.1em}\mathop{\leqslant}_{\LC\RC}\nolimits\hspace{0.1em}}

\def\lelrsn{\hspace{0.1em}\mathop{\leqslant}_{\LC\RC}^{\SG_{l,n-l}}
\nolimits\hspace{0.1em}}
\def\lelsn{\hspace{0.1em}\mathop{\leqslant}_\LC^{\SG_{l,n-l}}
\nolimits\hspace{0.1em}}
\def\lelrsl{\hspace{0.1em}\mathop{\leqslant}_{\LC\RC}^{\SG_l}
\nolimits\hspace{0.1em}}
\def\lelrwl{\hspace{0.1em}\mathop{\leqslant}_{\LC\RC}^{W_l}
\nolimits\hspace{0.1em}}

\def\lelo{\hspace{0.1em}\mathop{\leqslant}_\LC^\circ\nolimits\hspace{0.1em}}

\def\lero{\hspace{0.1em}\mathop{\leqslant}_\RC^\circ\nolimits\hspace{0.1em}}

\def\lelro{\hspace{0.1em}\mathop{\leqslant}_{\LC\RC}^\circ\nolimits\hspace{0.1em}}
\def\lepointo{\hspace{0.1em}\mathop{\leqslant}_?^\circ\nolimits\hspace{0.1em}}
\def\lepoint{\hspace{0.1em}\mathop{\leqslant}_?\nolimits\hspace{0.1em}}

\begin{document}

\baselineskip=16pt

\title{Two-sided cells in type ${\boldsymbol{B}}$ (asymptotic case)}

\author{C\'edric Bonnaf\'e}
\address{\noindent 
Labo. de Math. de Besan\c{c}on (CNRS: UMR 6623), 
Universit\'e de Franche-Comt\'e, 16 Route de Gray, 25030 Besan\c{c}on
Cedex, France} 

\makeatletter
\email{bonnafe@math.univ-fcomte.fr}

\makeatother


\date{\today}

\begin{abstract} 
We compute two-sided cells of Weyl groups of type $B$ for 
the ``asymptotic'' choice of parameters. We also obtain some 
partial results concerning Lusztig's conjectures in this 
particular case. 
\end{abstract}

\maketitle

\pagestyle{myheadings}

\markboth{\sc C. Bonnaf\'e}{\sc Two-sided cells in type $B$}


Let $W_n$ be a Weyl group of type $B_n$. The present paper is a continuation 
of the work done by L. Iancu and the 
author \cite{lacriced} concerning Kazhdan-Lusztig theory of $W_n$ 
for the asymptotic choice of parameters \cite[\SEC 6]{lacriced}. 
To each element $w \in W_n$ is associated 
a pair of standard bi-tableaux $(P(w),Q(w))$ 
(see \cite{okada} or \cite[\SEC 3]{lacriced}): 
this can be viewed as a Robinson-Schensted type correspondence. 
Our main result \cite[Theorem 7.7]{lacriced} 
was the complete determination of the left cells: two 
elements $w$ and $w'$ are in the same left cell if and only if 
$Q(w)=Q(w')$. For the corresponding result 
for the symmetric group, see \cite{KaLu} and \cite{ariki}. 
We have also computed 
the character afforded by a left cell representation 
\cite[Proposition 7.11]{lacriced} (this character is irreducible). 

In this paper, we are concerned with the computation of the 
two-sided cells. 
Let us state the result here. If $w \in W_n$, write 
$Q(w)=(Q^+(w),Q^-(w))$ and denote by 
$\l^+(w)$ and $\l^-(w)$ the shape of $Q^+(w)$ and $Q^-(w)$ respectively. 
Note that $(\l^+(w),\l^-(w))$ is a bipartition of $n$. 

\medskip

\noindent{\bf Theorem  (see \ref{cellules partitions}).} 
{\it For the choice of parameters as in \cite[\SEC 6]{lacriced}, 
two elements $w$ and $w'$ are in the same two-sided cell if and 
only if $(\l^+(w),\l^-(w))=(\l^+(w'),\l^-(w'))$.}

\medskip

Lusztig \cite[Chapter 14]{lusztig} has proposed fifteen 
conjectures on Kazhdan-Lusztig theory of Hecke algebras 
with unequal parameters. In the asymptotic case, Geck and Iancu \cite{lacrimeinolf} 
use some of our results, namely some informations on the preorder 
$\leqslant_{\LC\RC}$ (see Theorem \ref{ordre lr} and Proposition 
\ref{decroissante}), to compute the function $\ab$ 
and to prove Lusztig's conjectures 
$P_i$, for $i \in \{1,2,3,4,5,6,7,8,11,12,13,14\}$. 
On the other hand, Geck \cite{geck 2} has shown that Lusztig's conjectures 
$P_9$ and $P_{10}$ hold. More precisely, he proved that the Kazhdan-Lusztig 
basis is {\it cellular} (in the sense of \cite{lehrer}). 
He also proved a slightly weaker version 
of $P_{15}$ (but his version is sufficient for constructing the 
homomorphism from the Hecke algebra to the asymptotic algebra $J$).

\medskip

The present paper is organized as follows. In Section \ref{general}, 
we study some consequences of Lusztig's conjectures on 
the multiplication by $T_{w_0}$, where $w_0$ is the longest 
element of a finite Weyl group. From Section \ref{asymptotic} to the 
end of the paper, we assume that the Weyl group is of type $B_n$ 
and that the choice of parameters is done as in \cite[\SEC 6]{lacriced}. 
In Section \ref{asymptotic}, we establish some preliminary results 
concerning the Kazhdan-Lusztig basis. In Section \ref{bi}, 
we prove the above theorem by introducing a new basis of the 
Hecke algebra: this was inspired by the work of Geck on 
the induction of Kazhdan-Lusztig cells 
\cite{geck}. Section \ref{conj} contains some results related to 
Lusztig's conjectures. In Section \ref{specialization}, 
we determine which specializations of the parameters 
preserve the Kazhdan-Lusztig basis.

\medskip

\section{Generalities\label{general}}~

\medskip

\subsection{Notation} 
We slightly modify the notation used in \cite[\SEC 5]{lacriced}. 
Let $(W,S)$ be a Coxeter group with $|S| < \infty$. We denote by 
$\ell : W \to \NM=\{0,1,2,\dots\}$ the length function relative to $S$. 
If $W$ is finite, $w_0$ denotes its longest element. Let $\le$ denote 
the Bruhat ordering on $W$. If $I \incl S$, we denote by $W_I$ 
the standard parabolic subgroup of $W$ generated by $I$.

Let $\G$ be a totally ordered abelian group which will be denoted 
additively. The order on $\G$ will be denoted by $\le$. 
If $\g_0 \in \G$, we set 
$$\G_{< \g_0} = \{\g \in \G~|~\g < \g_0\},\qquad 
\G_{\le \g_0} = \{\g \in \G~|~\g \le \g_0\},$$
$$\G_{> \g_0} = \{\g \in \G~|~\g > \g_0\}\qquad \text{and}\qquad 
\G_{\ge \g_0} = \{\g \in \G~|~\g \ge \g_0\}.$$
Let $A$ be the group algebra of $\G$ over $\ZM$.  
It will be denoted exponentially: as a $\ZM$-module, it is free 
with basis $(v^\g)_{\g \in \G}$ and the multiplication rule 
is given by $v^\g v^{\g'}=v^{\g+\g'}$ for all $\g$, $\g' \in \G$. 
If $a \in A$, we denote by $a_\g$ the coefficient of $a$ on $v^\g$, so 
that $a=\sum_{\g \in \G} a_\g v^\g$. If $a \not= 0$, we define the {\it degree} 
and the {\it valuation} of $a$ (which we denote respectively by 
$\deg a$ and $\val a$) as the elements of $\G$ equal to 
$$\deg a=\max\{\g~|~a_\g \not= 0\}$$
$$\val a=\min\{\g~|~a_\g \not= 0\}.\leqno{\text{and}}$$
By convention, we set $\deg 0 = - \infty$ and $\val 0 = + \infty$. So 
$\deg : A \to \G \cup \{-\infty\}$ and $\val : A \to \G \cup \{+\infty\}$ 
satisfy $\deg ab = \deg a + \deg b$ and $\val ab=\val a + \val b$ 
for all $a$, $b \in A$. We denote by $A \to A$, $a \mapsto \bar{a}$ 
the automorphism of $A$ induced by the automorphism of $\G$ 
sending $\g$ to $-\g$. Note that $\deg a = -\val \aba$. 
If $\g_0 \in \G$, we set 
$$A_{< \g_0} = \Oplus_{\g < \g_0} \ZM v^\g,\qquad 
A_{\le \g_0} = \Oplus_{\g \le \g_0} \ZM v^\g,$$
$$A_{> \g_0} = \Oplus_{\g > \g_0} \ZM v^\g\qquad \text{and}\qquad
A_{\ge \g_0} = \Oplus_{\g \ge \g_0} \ZM v^\g.$$

We fix a {\it weight function} $L : W \to \G$, that is 
a function satisfying $L(ww')=L(w)+L(w')$ whenever 
$\ell(ww')=\ell(w)+\ell(w')$. We also assume that $L(s) > 0$ 
for every $s \in S$. We denote by $\HC=\HC(W,S,L)$ the 
{\it Hecke algebra} of $W$ associated to the weight function $L$.  
It is the associative $A$-algebra with $A$-basis $(T_w)_{w \in W}$ indexed 
by $W$ and whose multiplication is determined by the following two conditions:
$$\begin{array}{ccll}
(\text{a}) & \qquad & T_w T_{w'} =T_{ww'} & 
\text{if } \ell(ww')=\ell(w)+\ell(w')\qquad \\
(\text{b}) & \qquad & T_s^2 = 1 + (v^{L(s)}-v^{-L(s)}) T_s \quad& 
\text{if } s \in S. \\
\end{array}$$
It is easily seen from the above relations that $(T_s)_{s \in S}$ generates  
the $A$-algebra $\HC$ and that $T_w$ is invertible for every $w \in W$. 
If $h=\sum_{w \in W} a_w T_w \in \HC$, we set 
$\hov=\sum_{w \in W} \overline{a}_w T_{w^{-1}}^{-1}$. Then the 
map $\HC \to \HC$, $h \mapsto \hov$ is a semi-linear involutive 
automorphism of $\HC$. If $I \incl S$, we denote by $\HC(W_I)$ 
the sub-$A$-algebra of $\HC$ generated by $(T_s)_{s \in I}$. 

Let $w \in W$. By \cite[Theorem 5.2]{lusztig}, there exists a unique 
element $C_w \in \HC$ such that 
$$\begin{array}{ccl}
(\text{a}) & \qquad & C_w=\overline{C}_w \\
(\text{b}) & \qquad & C_w \in T_w + 
\Bigl(\DS{\mathop{\oplus}_{y \in W}} A_{< 0} T_y\Bigr).\qquad\qquad\\
\end{array}$$
Write $C_w=\sum_{y \in W} p_{y,w}^* T_y$ with $p_{y,w}^* \in A$. 
Then \cite[5.3]{lusztig}
$$\begin{array}{ll}
p_{w,w}^*=1& \\
p_{y,w}^* = 0 &\qquad \text{if }y \notle w.
\end{array}$$
In particular, $(C_w)_{w \in W}$ is an $A$-basis of $\HC$: it is called 
the {\it Kazhdan-Lusztig basis} of $\HC$. 
Write now $p_{y,w}=v^{L(w)-L(y)} p_{y,w}^*$. Then 
$$p_{y,w} \in A_{\ge 0}$$
and the coefficient of $p_{y,w}$ on $v^0$ is equal to $1$ 
(see \cite[Proposition 5.4 (a)]{lusztig}. 

We define the relations $\lel$, $\ler$, $\lelr$, $\sim_\LC$, $\sim_\RC$ 
and $\sim_{\LC\RC}$ as in \cite[\SEC 8]{lusztig}. 

\bigskip

\subsection{The function $\ab$} 
Let $x$, $y \in W$. Write 
$$C_x C_y =\sum_{z \in W} h_{x,y,z} C_z,$$
where $h_{x,y,z} \in A$ for $z \in W$. Of course, we have 
\equat\label{bar h}
\overline{h_{x,y,z}}=h_{x,y,z}.
\endequat

The following lemma is well-known 
\cite[Lemma 10.4 (c) and formulas 13.1 (a) et (b)]{lusztig}.

\medskip

\begin{lem}\label{borne degre}
Let $x$, $y$ and $z$ be three elements of $W$. Then 
$$\deg h_{x,y,z} \le \min(L(x),L(y)).$$
\end{lem}

\bigskip

\noindent{\bf Conjecture ${\boldsymbol{P}}_{\boldsymbol{0}}$ (Lusztig): } 
{\it There exists $N \in \G$ such that $\deg h_{x,y,z} \le N$ for all 
$x$, $y$ and $z$ in $W$.}

\bigskip

If $W$ is finite, then $W$ satisfies obviously $P_0$. If $W$ is an 
affine Weyl group, then it also satisfies $P_0$ \cite[7.2]{lusztig affine}. 
From now on, we assume that $W$ satisfies $P_0$, so that the next definition 
is valid. If $z \in W$, we set 
$$\ab(z)=\max_{x,y \in W} \deg h_{x,y,z}.$$
Since $h_{1,z,z}=1$, we have $\ab(z) \in \G_{\ge 0}$. If necessary, 
we will write $\ab_W(z)$ for $\ab(z)$. 
We denote by $\g_{x,y,z^{-1}} \in \ZM$ the coefficient of $v^{\ab(z)}$ in 
$h_{x,y,z}$. The next proposition shows how the function $\ab$ can be calculated 
by using different bases. 

\begin{prop}\label{changement de base}
Let $(X_w)_{w \in W}$ and $(Y_w)_{w \in W}$ be two families of elements of 
$\HC$ such that, for every $w \in W$, $X_w-T_w$ and $Y_w-T_w$ belong to 
$\mathop{\oplus}_{y < w} A_{<0} T_y$. For all $x$ and $y$ 
in $W$, write 
$$X_x Y_y = \sum_{z \in W} \x_{x,y,z} C_z.$$
Then, if $x$, $y$, $z \in W$, we have:
\begin{itemize}
\itemth{a} $\deg \x_{x,y,z} \le \min\{L(x),L(y)\}$.

\itemth{b} $\x_{x,y,z} \in \g_{x,y,z^{-1}} v^{\ab(z)} + A_{<\ab(z)}$. 
\end{itemize}
\noindent In particular, 
$$\ab(z)=\max_{x,y \in W} \deg \x_{x,y,z}.$$
\end{prop}

\proof Clear.\fin

\bigskip

\subsection{Lusztig's conjectures} 
Let $\t : \HC \to A$ be the $A$-linear map such that 
$\t(T_w)=\d_{1,w}$ if $w \in W$. 
It is the canonical symmetrizing form on $\HC$ (recall that 
$\t(T_x T_y)=\d_{xy,1}$). If $z \in W$, let 
$$\D(z)=-\deg p_{1,z}^*=-\deg \t(C_z).$$
Let $n_z$ be the coefficient of $p_{1,z}^*$ on $v^{-\D(z)}$. Finally, let 
$$\DC=\{z \in W~|~\ab(z) = \D(z)\}.$$

\noindent{\bf Conjectures (Lusztig): } {\it With the above notation, we have:

\begin{itemize}
\item[{\bf ${\boldsymbol{P}}_{\boldsymbol{1}}$.}] If $z \in W$, then 
$\ab(z) \le \D(z)$.

\item[{\bf ${\boldsymbol{P}}_{\boldsymbol{2}}$.}] If $d \in \DC$ and if $x$, 
$y \in W$ satisfy $\g_{x,y,d} \not= 0$, then $x=y^{-1}$.

\item[{\bf ${\boldsymbol{P}}_{\boldsymbol{3}}$.}] If $y \in W$, then 
there exists a unique $d \in \DC$ such that $\g_{y^{-1},y,d} \not= 0$.

\item[{\bf ${\boldsymbol{P}}_{\boldsymbol{4}}$.}] If $z' \lelr z$, then 
$\ab(z) \le \ab(z')$. Therefore, if $z \sim_{\LC\RC} z'$, then 
$\ab(z)=\ab(z')$. 

\item[{\bf ${\boldsymbol{P}}_{\boldsymbol{5}}$.}] If $d \in \DC$ and 
$y \in W$ satisfy $\g_{y^{-1},y,d} \not= 0$, then 
$\g_{y^{-1},y,d}=n_d=\pm 1$.

\item[{\bf ${\boldsymbol{P}}_{\boldsymbol{6}}$.}] If $d \in \DC$, then $d^2=1$.

\item[{\bf ${\boldsymbol{P}}_{\boldsymbol{7}}$.}] If $x$, $y$, $z \in W$, 
then $\g_{x,y,z}=\g_{y,z,x}$. 

\item[{\bf ${\boldsymbol{P}}_{\boldsymbol{8}}$.}] If $x$, $y$, $z \in W$ 
satisfy $\g_{x,y,z} \not= 0$, then $x \sim_\LC y^{-1}$, $y \sim_\LC z^{-1}$ and 
$z \sim_\LC x^{-1}$. 

\item[{\bf ${\boldsymbol{P}}_{\boldsymbol{9}}$.}] If $z' \lel z$ and 
$\ab(z')=\ab(z)$, then $z' \sim_\LC z$.

\item[{\bf ${\boldsymbol{P}}_{\boldsymbol{10}}$.}] If $z' \ler z$ and 
$\ab(z')=\ab(z)$, then $z' \sim_\RC z$.

\item[{\bf ${\boldsymbol{P}}_{\boldsymbol{11}}$.}] If $z' \lelr z$ and 
$\ab(z')=\ab(z)$, then $z' \sim_{\LC\RC} z$.

\item[{\bf ${\boldsymbol{P}}_{\boldsymbol{12}}$.}] If $I \incl S$ and 
$z \in W_I$, then $\ab_{W_I}(z)=\ab_W(z)$.

\item[{\bf ${\boldsymbol{P}}_{\boldsymbol{13}}$.}] Every left cell 
$\CC$ of $W$ contains a unique element $d \in \DC$. If 
$y \in \CC$, then $\g_{y^{-1},y,d} \not= 0$.

\item[{\bf ${\boldsymbol{P}}_{\boldsymbol{14}}$.}] If $z \in W$, then 
$z \sim_{\LC\RC} z^{-1}$. 

\item[{\bf ${\boldsymbol{P}}_{\boldsymbol{15}}$.}] If $x$, $x'$, $y$, $w \in W$ 
are such that $\ab(y)=\ab(w)$, then 
$$\sum_{y' \in W} h_{w,x',y'} \otimes_\ZM h_{x,y',y} = 
\sum_{y' \in W} h_{y',x',y} \otimes_\ZM h_{x,w,y'}$$
in $A \otimes_\ZM A$.
\end{itemize}}

\bigskip

Lusztig has shown that these conjectures hold if $W$ is a finite 
or affine Weyl group and  $L=\ell$ \cite[\SEC 15]{lusztig}, if 
$W$ is dihedral and $L$ is any weight function \cite[\SEC 17]{lusztig} 
and if $(W,L)$ is quasi-split \cite[\SEC 16]{lusztig}. 

\bigskip

\subsection{Lusztig's conjectures and multiplication by ${\boldsymbol{T_{w_0}}}$} 
We assume in this subsection that $W$ is finite. 
We are interested here in certain properties 
of the multiplication by $T_{w_0}^n$ for $n \in \ZM$. 
Some of them are partially known 
\cite[Lemma 1.11 and Remark 1.12]{lusztig symplectic}. 
If $y \in W$ and $n \in \ZM$, we set 
$$T_{w_0}^n C_y=\sum_{x \in W} \l^{(n)}_{x,y} C_x.$$
Note that $\l^{(n)}_{x,y}=0$ if $x \notlel y$. 

\medskip

\begin{prop}\label{multiplication tw0}
Assume that $W$ is finite and satisfies Lusztig's conjectures $P_1$, $P_4$ and 
$P_8$. Let $n \in \ZM$ and let $x$ and $y$ 
be two elements of $W$ such that $x \lel y$. Then: 

\begin{itemize}
\itemth{a}
If $n \ge 0$, then $\deg \l^{(n)}_{x,y} \le n (\ab(x)-\ab(w_0x))$. If 
moreover $x <_\LC y$, then $\deg \l^{(n)}_{x,y} < n (\ab(x)-\ab(w_0x))$.

\itemth{b} 
If $n \le 0$, then $\deg \l^{(n)}_{x,y} \le n (\ab(y)-\ab(w_0y))$. If 
moreover $x <_\LC y$, then $\deg \l^{(n)}_{x,y} < n (\ab(y)-\ab(w_0y))$.

\itemth{c} If $n$ is even and if $x \sim_\LC y$, then 
$\l^{(n)}_{x,y}=\d_{x,y} v^{n (\ab(x)-\ab(w_0x))}$. 
\end{itemize}
\end{prop}

\proof If $n=0$, then (a), (b) and (c) are easily checked. 
Let us now prove (a) and (b). By \cite[Proposition 11.4]{lusztig}, 
$$T_{w_0}=\sum_{u \in W} (-1)^{\ell(w_0u)} p_{1,w_0u}^* C_u.$$
Consequently, 
$$\l_{x,y}^{(1)}=\sum_{\substack{u \in W \\ x \ler u}} 
(-1)^{\ell(w_0u)} p_{1,w_0u}^* h_{u,y,x}.$$
But, by $P_1$, we have $\deg p_{1,w_0u}^* \le - \ab(w_0 u)$. 
If moreover $x \ler u$, then $w_0 u \ler w_0 x$ and so 
$-\ab(w_0 x) \ge -\ab(w_0 u)$ by $P_4$. Therefore 
$$\deg \l_{x,y}^{(1)} \le \ab(x)-\ab(w_0 x).$$
On the other hand, if $\deg \l_{x,y}^{(1)} = \ab(x)-\ab(w_0 x)$, 
then there exists $u \in W$ such that $x \ler u$ and $\deg h_{u,y,x} = \ab(x)$. 
So, by $P_8$, we get that $x \sim_\LC y$. This shows (a) for $n=1$. 

Now, let $\nu : \HC \to \HC$ denote the $A$-linear 
map such that 
$\nu(C_w)=v^{\ab(w_0w)-\ab(w)} C_w$ for all $w \in W$ and let 
$\mu : \HC \to \HC$, $h \mapsto T_{w_0} h$. Then, if $w \in W$, we have 
$$\nu\mu(C_y)= \sum_{u \lel y} v^{\ab(w_0u)-\ab(u)} \l_{u,y}^{(1)} C_u.$$
So, by the previous discussion, 
we have $v^{\ab(w_0u)-\ab(u)} \l_{u,y}^{(1)} \in A_{\le 0}$. 
Moreover, if $u <_\LC y$, then 
$v^{\ab(w_0u)-\ab(u)} \l_{u,y}^{(1)} \in A_{< 0}$ . 
On the other hand, $\det \mu=\pm 1$ and $\det \nu=1$. 
Therefore, if we write 
$$\mu^{-1}\nu^{-1}(C_y)=\sum_{u \lel y} \b_{u,y} C_u,$$
then $\b_{u,y} \in A_{\le 0}$ and, if $u <_\LC y$, then 
$\b_{u,y} \in A_{< 0}$. Finally, 
\eqna
T_{w_0}^{-1} C_y&=&\mu^{-1}(C_y)\\
   &=&\mu^{-1}\nu^{-1}\nu(C_y)\\
   &=&v^{\ab(w_0y)-\ab(y)} \mu^{-1}\nu^{-1}(C_y) \\
   &=&\sum_{u \lel y} v^{\ab(w_0y)-\ab(y)} \b_{u,y} C_u.
\endeqna
In other words, $\l_{x,y}^{(-1)}=v^{\ab(w_0y)-\ab(y)} \b_{x,y}$. 
This shows that (b) holds if $n=-1$. 
An elementary induction argument using $P_4$ shows that (a) and (b) 
hold in full generality. 

\medskip

Let us now prove (c). Let $K$ be the field of fraction of $A$. Let 
$C$ be a left cell of $W$ and let $c \in C$. We set 
$$\HC^{\lel C}=\mathop{\oplus}_{w \lel c} A C_w 
\qquad\text{and}\qquad 
\HC^{<_\LC C}=\mathop{\oplus}_{w <_\LC c} A C_w.$$
Then $\HC^{\lel C}$ and $\HC^{<_\LC C}$ are left ideals of 
$\HC$. The algebra $K\HC= K \otimes_A \HC$ being semisimple, 
there exists a left ideal $I_C$ of $K\HC$ 
such that $K\HC^{\lel C} = K\HC^{<_\LC C} \oplus I_C$. 

We need to prove that, for all $h \in I_C$, 
$$T_{w_0}^n h = v^{n (\ab(c)-\ab(w_0c))} h.$$
For this, we may, and we will, assume that $n > 0$. 
Let $V_1^C$, $V_2^C$,\dots, $V_{n_C}^C$ be irreducible 
sub-$K\HC \otimes_A K$-modules of $I_C$ such that 
$$I_C=V_1^C \oplus \dots \oplus V_{n_C}^C.$$
Let $j \in \{1,2,\dots,n_C\}$. 
Since $T_{w_0}^n$ is central and invertible in $\HC$, there exists 
$\e \in \{1,-1\}$ and $i_j^C \in \G$ such that 
$$T_{w_0}^n h = \e v^{i_j^C} h$$
for every $h \in V_j^C$. By specializing $v^\g \mapsto 1$, we get 
that $\e=1$. Moreover, by (a) and (b), $i_j^C \le n (\ab(c)-\ab(w_0c))$. 
On the other hand, since $\det \mu = \pm 1$, we have $\det \mu^n = 1$. 
But $\det \mu^n = v^r$, where
\eqna
r&=&\sum_{C \in \LC\CC(W)} \sum_{j=1}^{n_C} i_j^C \dim V_j^C \\
&\le & n \sum_{C \in \LC\CC(W)} 
(\ab(C)-\ab(w_0C))\sum_{j=1}^{n_C} \dim V_j^C \\
&=& n \sum_{C \in \LC\CC(W)} (\ab(C)-\ab(w_0C)) |C| \\
&=& n \sum_{w \in W} (\ab(w)-\ab(w_0w)) \\
&=& 0.
\endeqna
Here, $\LC\CC(W)$ denotes the set of left cells in $W$ and, if $C \in \LC\CC(W)$, 
$\ab(C)$ denotes the value of $\ab$ on $C$ (according to $P_4$). 
The fact that $r=0$ forces 
the equality $i_j^C = n (\ab(C)-\ab(w_0C))$ for every left cell $C$ 
and every $j \in \{1,2,\dots,n_C\}$.\fin

\bigskip

\remark{?} 
Assume here that $w_0$ is central in $W$ and keep the notation 
of the proof of Proposition \ref{multiplication tw0} (c). Let 
$j \in \{1,2,\dots,n_C\}$. Then there exists 
$\e_j(C) \in \{1,-1\}$ et $e_j(C) \in \G$ such that 
$T_{w_0} h = \e_j(C) v^{e_j(C)} h$
for every $h \in V_j^C$. 

\smallskip

\noindent{\bf Question: } {\it Let $j$, $j' \in \{1,2,\dots,n_C\}$. 
Does $\e_j(C)=\e_{j'}(C)$~?}

\smallskip

A positive answer to this question would allow to generalize 
Proposition \ref{multiplication tw0} (c) to the case where 
$w_0^n$ is central.\finl

\medskip

\begin{coro}\label{tw0 cw}
Assume that $W$ is finite and satisfies Lusztig's conjectures 
$P_1$, $P_2$, $P_4$, $P_8$, $P_9$ and $P_{13}$. Let $w \in W$ 
and let $n \in \NM$. Then 
$\deg \t(T_{w_0}^{-n} C_w) \le -\ab(w)+n(\ab(w_0 w) - \ab(w))$. 
Moreover, $\deg \t(T_{w_0}^{-n} C_w) = -\ab(w)+n(\ab(w_0 w) - \ab(w))$ 
if and only if $w_0^n w^{-1} \in \DC$.
\end{coro}

\proof 
Assume first that $n$ is even. In particular, $w_0^n=1$. 
By Proposition \ref{multiplication tw0}, we have 
$$\t(T_{w_0}^{-n} C_w) = v^{n(\ab(w_0w)-\ab(w))} \t(C_w) + \sum_{x <_\LC w} 
\l_{x,w}^{(-n)} \t(C_x).$$
But, if $x <_\LC w$, then $\deg \t(C_x) = -\D(x) \le -\ab(x) \le -\ab(w)$ 
by $P_1$ and $P_4$. So, by Proposition \ref{multiplication tw0} (b), we 
have that $\deg \l_{x,w}^{(-n)} \t(C_x) < -\ab(w)+n(\ab(w_0 w) - \ab(w))$. 
Moreover, again by $P_1$, we have 
$\deg \t(C_w) =-\D(w) \le -\ab(w)$. 
This shows that 
$\deg \t(T_{w_0}^{-n} C_w) \le -\ab(w)+n(\ab(w_0 w) - \ab(w))$ 
and that equality holds if and only if $\D(w)=\ab(w)$, that is, if and only if 
$w \in \DC$, as desired. 

\medskip

Assume now that $n=2k+1$ for some natural number $k$. 
Recall that, by \cite[Proposition 11.4]{lusztig}, 
$T_{w_0}=\sum_{u \in W} (-1)^{\ell(w_0u)} p_{1,w_0u}^* C_u$. Therefore
\eqna
T_{w_0}^{-n} C_w&=&
\sum_{u \in W} (-1)^{\ell(w_0u)} p_{1,w_0u}^* T_{w_0}^{-n-1} C_u C_w\\
&=& \sum_{\substack{u,x \in W \\ x \lel w~\text{and}~x \ler u}} 
(-1)^{\ell(w_0u)} p_{1,w_0u}^* 
h_{u,w,x} T_{w_0}^{-n-1} C_x.
\endeqna
This implies that 
$$\t(T_{w_0}^{-n} C_w)=
\sum_{\substack{u,x \in W \\ x \lel w~\text{and}~x \ler u}} 
(-1)^{\ell(w_0u)} p_{1,w_0u}^* h_{u,w,x} \t(T_{w_0}^{-n-1}C_x),$$
$$\deg \t(T_{w_0}^{-n} C_w) \le 
\max_{\substack{u,x \in W \\ x \lel w~\text{and}~x \ler u}} 
\deg(p_{1,w_0u}^* h_{u,w,x} \t(T_{w_0}^{-n-1}C_x)).\leqno{\text{so}}$$
Let $u$ and $x$ be two elements of $W$ such that $x \lel w$ and 
$x \ler u$. Since $n+1$ is even and by the previous discussion, we have 
$$\deg \t(T_{w_0}^{-n-1}C_x) \le -\ab(x)+(n+1)(\ab(w_0x)-\ab(x)).$$
By $P_1$ and $P_4$, $\deg p_{1,w_0u}^* \le -\ab(w_0u) \le -\ab(w_0x)$. 
Moreover, $\deg h_{u,w,x} \le \ab(x)$. Consequently, 
\eqna
\deg(p_{1,w_0u}^* h_{u,w,x} \t(T_{w_0}^{-n-1}C_x)) &\le& 
-\ab(x)+n(\ab(w_0x)-\ab(x)) \\
&\le& -\ab(w)+n(\ab(w_0w)-\ab(w)).
\endeqna
Moreover, equality holds if and only if $w_0 u \in \DC$, 
$\deg h_{u,w,x}=\ab(x)=\ab(w)$ and $x \in \DC$. 

We first deduce that 
$$\deg \t(T_{w_0}^{-n} C_w) \le -\ab(w)+n(\ab(w_0w)-\ab(w))$$ 
which is the first assertion of the proposition. 

\smallskip

Assume now that 
$\deg \t(T_{w_0}^{-n} C_w) = -\ab(w)+n(\ab(w_0w)-\ab(w))$. 
Then there exists $u$ and $x$ in $W$ such that $x \lel w$, 
$x \ler u$, $w_0 u \in \DC$, $\deg h_{u,w,x}=\ab(x)=\ab(w)$ and $x \in \DC$. 
Since $\deg h_{u,w,x}=\ab(x)$ and $x \in \DC$, we deduce from $P_2$ that 
$w=u^{-1}$, which shows that $w_0 w^{-1} \in \DC$. 

Conversely, assume that $w_0 w^{-1}$ belongs to $\DC$. To show that 
$\deg \t(T_{w_0}^{-n} C_w) = -\ab(w)+n(\ab(w_0w)-\ab(w))$, it is sufficient 
to show that there is a unique pair $(u,x)$ of elements of $W$ 
such that $x \lel w$, $x \ler u$, $w_0 u \in \DC$, 
$\deg h_{u,w,x}=\ab(x)=\ab(w)$ and $x \in \DC$. The existence follows 
from $P_{13}$ (take $u=w^{-1}$ and $x$ be the unique el\'ement of $\DC$ 
belonging to the left cell containing $w$). 
Let us now show unicity. Let $(u,x)$ be such a pair. 
Since $\deg h_{u,w,x}=\ab(x)$ and $x \in \DC$, we deduce from $P_2$ that 
$u=w^{-1}$. Moreover, since 
$\ab(x)=\ab(w)$ and $x \lel w$, we have $x \sim_\LC w$ by $P_9$. 
But, by $P_{13}$, $x$ is the unique element of $\DC$ 
belonging to the left cell containing $w$.\fin

\bigskip

\section{Preliminaries on type $B$ (asymptotic case)\label{asymptotic}}

From now on, we are working under the following hypothesis. 

\smallskip

\begin{quotation}
\noindent{\bf Hypothesis and notation:} {\it We assume now that 
$W=W_n$ is of type $B_n$, $n \ge 1$. We write $S=S_n=\{t,s_1,\dots,s_{n-1}\}$ 
as in \cite[\SEC 2.1]{lacriced}: the Dynkin diagram of $W_n$ is given by 
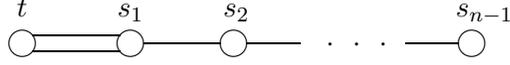
\begin{center}
\begin{picture}(220,30)
\put( 40, 10){\circle{10}}
\put( 44,  7){\line(1,0){33}}
\put( 44, 13){\line(1,0){33}}
\put( 81, 10){\circle{10}}
\put( 86, 10){\line(1,0){29}}
\put(120, 10){\circle{10}}
\put(125, 10){\line(1,0){20}}
\put(155,  7){$\cdot$}
\put(165,  7){$\cdot$}
\put(175,  7){$\cdot$}
\put(185, 10){\line(1,0){20}}
\put(210, 10){\circle{10}}
\put( 38, 20){$t$}
\put( 76, 20){$s_1$}
\put(116, 20){$s_2$}
\put(204, 20){$s_{n{-}1}$}
\end{picture}
\end{center}
We also assume that $\G=\ZM^2$ and that $\G$ is ordered lexicographically:
$$(a,b) \le (a',b') \Longleftrightarrow 
a < a' \text{ or } (a=a' \text{ and } b \le b').$$
We set $V=v^{(1,0)}$ and $v=v^{(0,1)}$ so that $A=\ZM[V,V^{-1},v,v^{-1}]$ 
is the Laurent polynomial ring in two algebraically independent indeterminates 
$V$ and $v$. If $w \in W_n$, we denote by $\ell_t(w)$ the number of occurences 
of $t$ in a reduced expression of $w$. We set $\ell_s(w)=\ell(w)-\ell_t(w)$. 
Then $\ell_s$ and $\ell_t$ are weight functions and we assume that
$L=L_n : W_n \to \G$, $w \mapsto (\ell_t(w),\ell_s(w))$. 
So $\HC=\HC_n=\HC(W_n,S_n,L_n)$. 
We denote by $\SG_n$ the subgroup of $W$ generated by $\{s_1,\dots,s_{n-1}\}$: 
it is isomorphic to the symmetric group of degree $n$.}
\end{quotation}

\medskip

We now recall some notation from \cite[\SEC 2.1 and 4.1]{lacriced}. 
Let $r_1=t_1=t$ and, if $1 \le i \le n-1$, let $r_{i+1} = s_i r_i$ 
and $t_{i+1} = s_i t_i s_i$. If $0 \le l \le n$, let $a_l = r_1 r_2 \dots r_l$. 
We denote by $\SG_l$, $W_l$, $\SG_{l,n-l}$, $W_{l,n-l}$ the standard 
parabolic subgroups of $W_n$ generated by $\{s_1,s_2,\dots,s_{l-1}\}$, 
$\{t,s_1,s_2,\dots, s_{l-1}\}$, $S_n \setminus \{t,s_l\}$ and 
$S_n \setminus \{s_l\}$ respectively. The longest element 
of $\SG_l$ is denoted by $\s_l$. Let 
$$Y_{l,n-l}=\{a \in \SG_n~|~\forall~\s \in \SG_{l,n-l},~\ell(a\s) \ge \ell(\s)\}.$$
If $w \in W_n$ is such that $\ell_t(w)=l$, then \cite[\SEC 4.6]{lacriced} 
there exist unique $a_w$, $b_w \in Y_{l,n-l}$, 
$\s_w \in \SG_{l,n-l}$ such that $w=a_w a_l \s_w b_w^{-1}$. Recall that 
$\ell(w)=\ell(a_w)+\ell(a_l)+\ell(\s_w)+\ell(b_w)$. 

\bigskip

\subsection{Some submodules of $\HCB$\label{soussection ideaux}} 
If $l$ is a natural number such that $0 \le l \le n$, we set 
$$\TC_l=\mathop{\oplus}_{\substack{w \in W_n \\ \ell_t(w) = l}} A T_w,\qquad
\TC_{\le l}=\mathop{\oplus}_{\substack{w \in W_n \\ \ell_t(w) \le l}} A T_w,\qquad
\TC_{\ge l}=\mathop{\oplus}_{\substack{w \in W_n \\ \ell_t(w) \ge l}} A T_w,$$
$$\CC_l=\mathop{\oplus}_{\substack{w \in W_n \\ \ell_t(w) = l}} A C_w,\qquad
\CC_{\le l}=\mathop{\oplus}_{\substack{w \in W_n \\ \ell_t(w) \le l}} A C_w\quad
\text{and}\quad
\CC_{\ge l}=\mathop{\oplus}_{\substack{w \in W_n \\ \ell_t(w) \ge l}} A C_w.$$
Let $\Pi_?^T : \HC_n \to \TC_?$ and 
$\Pi_?^C : \HC_n \to \CC_?$ 
be the natural projections (for $? \in \{l,\le l, \ge l\}$). 

\medskip

\begin{prop}\label{ideaux}
Let $l$ be a natural number such that $0 \le l \le n$. Then:
\begin{itemize}
\itemth{a} $\TC_l$ and $\CC_l$ are 
sub-$\HC(\SG_n)$-modules-$\HC(\SG_n)$ of $\HC_n$. The maps 
$\Pi_l^T$ and $\Pi_l^C$ are morphisms of 
$\HC(\SG_n)$-modules-$\HC(\SG_n)$.

\itemth{b} $\TC_{\le l}=\CC_{\le l}$. 

\itemth{c} $\CC_{\ge l}$ is a two-sided ideal of $\HC_n$. 
\end{itemize}
\end{prop}

\proof 
(a) follows from \cite[Theorem 6.3 (b)]{lacriced}.
(b) is clear.
(c) follows from \cite[Corollary 6.7]{lacriced}.\fin

\bigskip

The next proposition is a useful characterization of the elements 
of the two-sided ideal $\CC_n$. 

\medskip

\begin{prop}\label{caracterisation}
Let $h \in \HC_n$. The following are equivalent:
\begin{itemize}
\itemth{1} $h \in \CC_n$.

\itemth{2} $\forall x \in \HC_n$, $(T_t-V) xh=0$.

\itemth{3} $\forall x \in \HC(\SG_n)$, $(T_t-V) xh=0$.
\end{itemize}
\end{prop}

\proof 
If $\ell_t(w)=n$, then $tw < w$ so $(T_t-V) C_w=0$ by 
\cite[Theorem 6.6 (b)]{lusztig}. Since $\CC_n$ is a two-sided ideal of $\HC_n$ 
(see Proposition \ref{ideaux} (c)), 
we get that (1) implies (2). It is also obvious that 
(2) implies (3). It remains to show that (3) implies (1).

\medskip

Let $I=\{h \in \HC_n~|~\forall x \in \HC(\SG_n)$, $(T_t-V) xh=0\}$. 
Then $I$ is clearly a sub-$\HC(\SG_n)$-module-$\HC_n$ of $\HC_n$.
We need to show that $I \incl \CC_n$. In other words, 
since $\CC_n \incl I$, we need to show that 
$I \cap \CC_{\le n-1}=0$. 
Let $I'=I \cap \CC_{\le n-1}=I \cap \TC_{\le n-1}$ 
(see Proposition \ref{ideaux} (b)). 
Let $X=\{w \in W_n~|~\t(I' T_{w^{-1}}) \not= 0\}$. Showing that $I'=0$ 
is equivalent to showing that $X=\vide$. 

Assume $X\not=\vide$. Let $w$ be an element of $X$ of maximal length 
and let $h$ be an element of $I'$ such that $\t(hT_{w^{-1}})\not= 0$. 
Since $h \in \TC_{\le n-1}$, we have $\ell_t(w) \le n-1$. 
Moreover, $T_t h = Vh$, so $\t(T_t h T_{w^{-1}})=\t(h T_{w^{-1}} T_t) \not = 0$. 
By the maximality of $\ell(w)$, we get that $tw < w$. 
So, there exists $s \in S_n$ such that $sw > w$ and $s \not= t$. 
Then 
\eqna
\t(T_sh T_{(sw)^{-1}}) &=& \t(h T_{(sw)^{-1}} T_s) \\
&=& \t(h T_{w^{-1}}) + (v - v^{-1}) \t(h T_{(sw)^{-1}}) \\
&\not=& 0,
\endeqna
the last inequality following from the maximality of $\ell(w)$ 
(which implies that $\t(h T_{(sw)^{-1}})=0$). 
But $T_s h\in I'$ and so $sw \in X$. This contradicts 
the maximality of $\ell(w)$.\fin
 
\bigskip

\subsection{Some results on the Kazhdan-Lusztig basis\label{formules}}
In this subsection, we study the elements of the Kazhdan-Lusztig basis 
of the form $C_{a_l \s}$ where $0 \le l \le n$ and $\s \in \SG_n$. 

\medskip

\begin{prop}\label{cal csigma}
Let $\s \in \SG_n$ and let $0 \le l \le n$. Then 
$C_{a_l} C_\s=C_{a_l \s}$ and $C_\s C_{a_l}=C_{\s a_l}$. 
\end{prop}

\proof
Let $C=C_{a_l}C_\s$. Then $\Cov=C$ and 
$$C-T_{a_l\s} =\sum_{w < a_l \s} \l_w T_w$$
with $\l_w \in A$ for $w < a_l \s$. To show that 
$C=C_{a_l\s}$, it is sufficient to show that $\l_w \in A_{< 0}$. 

But, $C_{a_l}=T_{a_l} + \sum_{x < a_l} V^{\ell_t(x)-l} \b_x T_x$ with  
$\b_x \in \ZM[v,v^{-1}]$ (see \cite[Theorem 6.3 (a)]{lacriced}). 
Hence, 
$$C=T_{a_l \s} + \sum_{\t < \s} p_{\t,\s}^* T_{a_l \t} + 
\sum_{x < a_l} V^{\ell_t(x)-l} \b_x T_x C_\s.$$
But, if $x < a_l$, then $\ell_t(x) < l$. This shows that 
$\l_w \in A_{<0}$ for every $w < a_l \s$.
This shows the first equality. The second one is obtained by a symmetric 
argument.\fin 

\bigskip

Proposition \ref{cal csigma} shows that it can be useful to compute 
in different ways the elements $C_{a_l}$ to be able to relate 
the Kazhdan-Lusztig basis of $\HC_n$ to the Kazhdan-Lusztig basis 
of $\HC(\SG_n)$. Following the work of Dipper-James-Murphy 
\cite{DiJa95}, Ariki-Koike \cite{ArKo} and 
Graham-Lehrer \cite[\SEC 5]{lehrer}, we set 
\eqna
P_l&=&(T_{t_1} + V^{-1})(T_{t_2}+V^{-1}) \dots (T_{t_l}+V^{-1}) \\
&=&\DS{\sum_{0 \le k \le l} V^{k-l} 
\bigl(\sum_{1 \le i_1 < \dots < i_k \le l} T_{t_{i_1} \dots t_{i_k}}\bigr).}
\endeqna

\medskip

\begin{lem}
$P_n$ is central in $\HC_n$. 
\end{lem}

\proof First, $P_n$ commutes with $T_t$ (indeed, $t t_i = t_i t > t_i$ for 
$1 \le i \le n$). By \cite[Lemma 3.3]{ArKo}, $P_n$ commutes with $T_{s_i}$ for 
$1 \le i \le n-1$. Since the notation and conventions are somewhat different, 
we recall here a brief proof. First, if $j \not\in \{i,i+1\}$, 
$s_i t_j = t_j s_i > t_j$ so $T_{s_i}$ commutes with $T_{t_j}$. Therefore, 
it is sufficient to show that $T_{s_i}$ commutes with 
$(T_{t_i}+V^{-1})(T_{t_{i+1}} + V^{-1})$. This follows from a straightforward 
computation using the fact that $s_i t_i > t_i$, that $t_{i+1} s_i < t_{i+1}$ and 
that $s_i t_i = t_{i+1} s_i$.\fin

\medskip

\begin{prop}\label{cal pl}
If $0 \le l \le n$, then $C_{a_l}=P_l T_{\s_l}^{-1}=T_{\s_l}^{-1} P_l$. 
\end{prop}

\proof 
The computation may be performed in the subalgebra of $\HC_n$ generated by 
$\{T_t,T_{s_1},\dots,T_{s_{l-1}}\}$ so we may, and we will, assume 
that $l=n$. First, we have $T_tP_n=VP_n$. Since $P_n$ is central 
in $\HC_n$, it follows from the characterization of $\CC_n$ 
given by Proposition \ref{caracterisation} that $P_n \in \CC_n$. 

Now, let $h=C_{a_n}-P_n T_{\s_n}^{-1}$. Then, by Proposition \ref{ideaux} (a), 
we have $h \in \CC_n$. Moreover, it is easily checked that 
$h \in \TC_{\le n-1}=\CC_{\le n-1}$. So $h=0$.\fin

\medskip

\begin{coro}\label{tau cal}
If $0 \le l \le n$ and $\s \in \SG_n$, then 
$\Pi_0^T(C_{a_l\s})=V^{-l}T_{\s_l}^{-1} C_\s$. 
In particular, $\t(C_{a_l\s})=V^{-l}\t(T_{\s_l}^{-1} C_\s)$. 
\end{coro}

\proof 
Since $\Pi_0^T$ is a morphism of right $\HC(\SG_n)$-modules 
(see Proposition \ref{ideaux} (a)) 
and since $C_{a_l\s}=P_l T_{\s_l}^{-1} C_\s$ (see Propositions 
\ref{cal csigma} and \ref{cal pl}), we have 
$\Pi_0^T(C_{a_l\s})=\Pi_0^T(P_l) T_{\s_l}^{-1} C_\s$. But, 
$\Pi_0^T(P_l)=V^{-l}$. This completes the proof of the corollary.\fin

\bigskip

\section{Two-sided cells\label{bi}}

\medskip

The aim of this section is to show that, if $x$ and $y$ are two elements 
of $W$ such that $\ell_t(x)=\ell_t(y)=l$, then $x \lelr y$ 
if and only if $\s_x \lelrsn \s_y$ (see Theorem \ref{ordre lr}). 
Here, $\lelrsn$ is the preorder $\lelr$ defined inside 
the parabolic subgroup $\SG_{l,n-l}$. 
For this, we adapt an argument of Geck \cite{geck} 
who was considering the preorder $\lel$. 

We start by defining an order relation $\infspe$ on $W$. 
Let $x$ and $y$ be two elements of $W$. Then $x \inferieur y$ 
if the following conditions are fulfilled:

\begin{quotation}
\noindent (1) $\ell_t(x)=\ell_t(y)$, 

\noindent (2) $x \le y$, 

\noindent (3) $a_x < a_y$ or $b_x < b_y$, 

\noindent (4) $\s_x \lelrsn \s_y$.
\end{quotation}

\smallskip

\noindent We write $x \infspe y$ if $x \inferieur y$ or $x=y$. 
If $y \in W_n$, we set 
$$\G_y=T_{a_y} C_{a_{\el_t(y)}} C_{\s_y} T_{b_y^{-1}}.$$

\medskip

\begin{lem}\label{inf 0}
Let $y \in W$. Then $\G_y \in T_y + \mathop{\oplus}_{x < y} A_{<0} T_x$.
\end{lem}

\proof First, recall that $\G_y$ is a linear combination of elements of the form 
$T_{a_y} T_ z T_{b_y^{-1}}$ with $z \le a_{\ell_t(y)} \s_y$, so it 
is a linear combination of elements of the form $T_x$ with $x \le y$. 

Let $l=\ell_t(y)$ and $\s=\s_y$. We have 
$$\G_y=T_{a_y} T_{a_l} T_\s T_{b_y^{-1}} + 
\bigl(\sum_{\t < \s} p_{\t,\s}^* T_{a_y} T_{a_l} T_\t T_{b_y^{-1}}\bigr) 
+ \bigl(\sum_{a < a_l,~\t \le \s} p_{a,a_l}^* p_{\t,\s}^* 
T_{a_y} T_a T_\t T_{b_y^{-1}}\bigr).$$
If $\t < \s$, then $T_{a_y} T_{a_l} T_\t T_{b_y^{-1}}=T_{a_y a_l \t b_y^{-1}}$ 
by \cite[\SEC 4.6]{lacriced}. 
On the other hand, if $a < a_l$, then $\ell_t(a) < l$ so 
$T_{a_y} T_a T_\t T_{b_y^{-1}}$ is a linear combination, with 
coefficients in $\ZM[v,v^{-1}]$ of elements $T_w$ with $\ell_t(w)=\ell_t(a) < l$ 
(because $a_y$, $\t$ and $b_y^{-1}$ are elements of $\SG_n$). 
Since $V^{l - \ell_t(a)} p_{a,a_l}^* \in \ZM[v,v^{-1}]$ by 
\cite[Theorem 6.3 (a)]{lacriced}, this proves the lemma.\fin

\medskip

\begin{lem}\label{geck}
If $y \in W_n$, then  
$$\overline{\G}_y=\G_y + \sum_{x \inferieur y} \r_{x,y} \G_x$$
where the $\r_{x,y}$'s belong to $\ZM[v,v^{-1}]$.
\end{lem}

\proof 
Let $l=\ell_t(y)$. Then 
$$T_{a_y^{-1}}^{-1}=T_{a_y} + 
\sum_{\substack{a \in Y_{l,n-l}\\ x \in \SG_{l,n-l} \\ ax < a_y}} R_{ax,a_y} T_a T_x.$$
Moreover, if $a \in Y_{l,n-l}$ and $x \in \SG_{l,n-l}$ are such that 
$ax < a_y$, then $a < a_y$ (see \cite[Lemma 9.10 (f)]{lusztig}. Thus,
$$\overline{\G}_y=\G_y + 
\sum_{\substack{a,b \in Y_{l,n-l} \\ x,x' \in \SG_{l,n-l} \\ ax < a_y 
~\text{or}~bx' < b_y}} R_{ax,a_y}R_{bx',b_y} 
T_a(T_x C_{a_l \s_y} T_{x^{\prime -1}}) T_{b^{-1}}.$$
The result now follows from Lemma \ref{inf 0}.\fin

\medskip

\begin{coro}\label{involution rho}
If $x \infspe y$, then 
$\DS{\sum_{x \infspe z \infspe y} \overline{\rho}_{x,z} \r_{z,y} = \d_{x,y}}$. 
\end{coro}

\proof 
This follows immediately from Lemma \ref{geck} and from the fact that 
$\HC \to \HC$, $h \mapsto \hov$ is an involution.\fin

\medskip

\begin{coro}\label{geck formule}
If $w \in W$, then 
$$C_w=\G_w + \sum_{y \inferieur w} \pi_{y,w}^* \G_y$$
where $\pi_{y,w}^* \in v^{-1}\ZM[v^{-1}] \incl A_{<0}$ if $y \inferieur w$.
\end{coro}

\proof By Corollary \ref{involution rho}, there exists a unique 
family $(\pi_{y,w}^*)_{y \inferieur w}$ of elements of $v^{-1}\ZM[v^{-1}]$ 
such that $\G_w + \sum_{y \inferieur w} \pi_{y,w}^* \G_y$ 
is stable under the involution $h \mapsto \hov$ of $\HC_n$ 
(see \cite[Page 214]{ducloux}: this contains a general setting for 
including the arguments in \cite[Proposition 2]{lusztig left} or in 
\cite[Proposition 3.3]{geck}). But, by Lemma \ref{inf 0}, we have 
$$\G_w + \sum_{y \inferieur w} \pi_{y,w}^* \G_y \in 
T_w + (\mathop{\oplus}_{y < w} A_{<0} T_y).$$
So $C_w=\G_w + \sum_{y \inferieur w} \pi_{y,w}^* \G_y$.\fin

\bigskip

We are now ready to prove the main theorem of this section. 

\medskip

\begin{theo}\label{ordre lr}
Let $x$ and $y$ be two elements of $W$ such that $\ell_t(x)=\ell_t(y)=l$. 
Then $x \lelr y$ if and only if $\s_x \lelrsn \s_y$. 
\end{theo}

\proof Assume first that 
$\s_x \lelrsn \s_y$. Decompose $\s_x=(\s_x',\s_x'')$ with $\s_x' \in \SG_l$ 
and $\s_x'' \in \SG_{n-l}$. Then $\s_x' \lelrsl \s_y'$ so 
$\s_l \s_y' \lelrsl \s_l \s_x'$ so $w_l \s_l \s_x' \lelrwl w_l \s_l \s_y'$. 
In other words, $a_l \s_x' \lelrwl a_l \s_y'$. Therefore, 
$a_l \s_x \lelr a_l \s_y$. But, by \cite[Theorem 7.7]{lacriced}, 
we have $x \sim_{\LC\RC} a_l \s_x$ and $y \sim_{\LC\RC} a_l \s_y$. So $x \lelr y$. 

\medskip

To show the converse statement, it is sufficient to show that 
$$I=\Bigl(\mathop{\oplus}_{\substack{u \in W_n \\ 
\ell_t(u)=l~\text{and}~\s_u \lelrsn \s_y}} 
AC_u\Bigr) \oplus \CC_{\ge l+1}$$
is a two-sided ideal. But, by Corollary \ref{geck formule}, 
we have 
$$I=\Bigl(\mathop{\oplus}_{\substack{u \in W_n \\ 
\ell_t(u)=l~\text{and}~\s_u \lelrsn \s_y}} 
A\G_u\Bigr) \oplus \CC_{\ge l+1}.$$

By symmetry, we only need to prove that $I$ is a left ideal. 
Let $h \in \HC_n$ and let $u \in W_n$ such that $\ell_t(u)=l$ 
and $\s_u \lelrsn \s_y$. We want to prove that $h \G_u \in I$. 
For simplification, let $a=a_u$, $b=b_u$, $\s=\s_u$. Let 
$$X_l = \{x \in W_n~|~\forall~w \in W_{l,n-l},~\ell(xw) \ge \ell(w)\}.$$
Then, by \cite[Proposition 3.3]{geck} and 
\cite[Lemma 7.3 and Corollary 7.4]{lacriced}, 
$$T_a C_{a_l \s} \in \mathop{\oplus}_{\substack{x \in X_l \\ 
\t \lelsn \s}} A C_{xa_l\t}.$$
Let $I'$ be the right-hand side of the previous formula. 
By \cite[Corollary 3.4]{geck}, $I'$ is a left ideal. Therefore, 
$h T_a C_{a_l \s} \in I'$. On the other hand, 
$$I' \subset \Bigl(\mathop{\oplus}_{\substack{x \in Y_{l,n-l} \\ 
\t \lelsn \s}} A C_{xa_l\t}\Bigr) \oplus \CC_{\ge l+1}.$$
Now, by Corollary \ref{geck formule}, we have 
$$I' \subset \Bigl(\mathop{\oplus}_{\substack{x \in Y_{l,n-l} \\ 
\t \lelsn \s}} A T_x C_{a_l\t}\Bigr) \oplus \CC_{\ge l+1}.$$
Therefore, $h \G_u \in I' T_{b^{-1}} \subset I$, as desired.\fin

\medskip

\begin{coro}\label{cellules bilateres}
Let $x$ and $y$ be two elements of $W_n$. Then $x \sim_{\LC\RC} y$ 
if and only if $\ell_t(x)=\ell_t(y)$ ($=l$) and 
$\s_x \sim_{\LC\RC}^{\SG_{l,n-l}} \s_y$.
\end{coro}

\bigskip

\remark{partitions} 
We associate to each element $w \in W_n$ a pair $(P(w),Q(w))$ 
of standard bi-tableaux as in \cite[\SEC 3]{lacriced}. Let $l=\ell_t(w)$. 
Write $Q(w)=(Q^+(w),Q^-(w))$ and denote by $\l^?(w)$ the shape of 
$Q^?(w)$ for $? \in \{+,-\}$. The map $w \mapsto (P(w),Q(w))$ 
is a generalization of the Robinson-Schensted correspondence 
(see \cite[Theorem 3.3]{okada} or \cite[Theoreme 3.3]{lacriced}). 
Then $\l^+(w)$ is a partition of $n-l$ 
and $\l^-(w)$ is a partition of $l$, so that $\l(w)=(\l^+(w),\l^-(w))$ 
is a bipartition of $n$. If we write $\s_w=\s_w^- \times \s_w^+$ with 
$\s_w^- \in \SG_l$ and $\s_w^+ \in \SG_{n-l}$, note that $\l^+(w)$ is 
the shape of the standard tableau associated to $\s_w^+$ 
by the classical Robinson-Schensted correspondence while 
$\l^-(w)^*$ (the partition conjugate to $\l^-(w)$) 
is the shape of the standard tableau associated 
to $\s_w^-$. Let $\trianglelefteq$ 
denote the {\it dominance order} on partitions: if 
$\a=(\a_1 \ge \a_2 \ge \dots)$ and $\b=(\b_1 \ge \b_2 \ge \dots)$ 
are two partitions of the same natural number, we write $\a \trianglelefteq \b$
if 
$$\sum_{j=1}^i \a_j \le \sum_{j=1}^i \b_j$$ 
for every $i \ge 1$. 
Now, let $x$ and $y$ be two elements of $W_n$. 
If $\ell_t(x)=\ell_t(y)$, then Theorem \ref{ordre lr} is equivalent to:
\equat\label{ordre partitions}
\text{\it $x \lelr y$ if and only if $\l^+(x) \trianglelefteq \l^+(y)$ 
and $\l^-(y) \trianglelefteq \l^-(x)$.}
\endequat
This follows from \cite[3.2]{LuXi} and \cite[2.13.1]{DPS} (see also 
\cite[Exercise 5.6]{ourbuch}). 
Then, for general $x$ and $y$, Corollary 
\ref{cellules bilateres} is equivalent to:
\equat\label{cellules partitions}
\text{\it $x \sim_{\LC\RC} y$ if and only if $\l(x)=\l(y)$.}
\endequat

\bigskip

\section{Around Lusztig's conjectures\label{conj}}

\medskip

In this section, we prove some results which are related to Lusztig's 
conjectures. If $\s \in \SG_n$, we denote by $\ab_\SG(\s)$ the function $\ab$ 
evaluated on $\s$ but computed in $\SG_n$. It is given by the following formula. 
Let $\l=(\l_1 \ge \l_2 \ge \dots )$ be the shape of the left 
cell of $\s$. Then
$$\ab_\SG(\s)=\sum_{i \ge 1} (i-1) \l_i.$$
We denote by $\ab_\l$ the right-hand side of the previous formula. 
If $z \in W$, we set
$$\alpb(z)=(\ell_t(z),2 \ab_\SG(\s_z)-\ab_\SG(\s_{\ell_t(z)}\s_z)) \in \NM^2.$$
In terms of partitions 
(using the notation introduced in Remark \ref{partitions}), we have
$$\alpb(z)=(|\l^-(z)|, \ab_{\l^+(z)} + 2\ab_{\l^-(z)^*} - \ab_{\l^-(z)}).$$
We now study some properties of the function $\alpb$. 

\medskip

\remark{a=a} 
Geck and Iancu \cite{lacrimeinolf} have proved, using the result 
of this section (and especially Proposition 4.2), that $\alpb=\ab$. 
They have deduced, using the notion of {\it orthogonal representations}, 
that Lusztig's conjectures $P_i$ hold for $i \in \{1,2,3,4,5,6,7,8,11,12,13,14\}$. 
After that, Geck \cite{geck 2} proved $P_9$ and $P_{10}$.\finl

\medskip

The first proposition shows that $\alpb$ is decreasing with respect 
to $\lelr$ (compare with Lusztig's conjecture $P_4$).

\medskip

\begin{prop}\label{decroissante}
Let $z$ and $z'$ be two elements of $W$. Then:
\begin{itemize}
\itemth{a} If $z \lelr z'$, then $\alpb(z') \le \alpb(z)$. 

\itemth{b} If $z \lelr z'$ and $\alpb(z)=\alpb(z')$ then $z \sim_{\LC\RC} z'$. 
\end{itemize}
\end{prop}

\proof Since $z \lelr z'$, we have $\ell_t(z) \ge \ell_z(z')$ 
by \cite[Corollary 6.7]{lacriced}. Therefore, if $\ell_t(z) > \ell_t(z')$, 
then $\alpb(z) > \alpb(z')$ and $z \not\sim_{\LC\RC} z'$. This proves 
(a) and (b) in this case. 

So, assume that $\ell_t(z)=\ell_t(z')=l$. Then, by Theorem 
\ref{ordre lr}, we have $\s_z \lelrsn \s_{z'}$. Write $\s_z=(\s,\t)$ 
and $\s_{z'}=(\s',\t')$ where $\s$, $\s' \in \SG_l$ 
and $\t$, $\t' \in \SG_{n-l}$. Then 
$$\alpb(z)=(l, 2\ab_\SG(\s) - \ab_\SG(\s_l \s) + \ab_\SG(\t))$$ 
$$\alpb(z')=(l, 2\ab_\SG(\s') - \ab_\SG(\s_l \s') + \ab_\SG(\t')).
\leqno{\text{and}}$$ 
But $\s \hspace{0.05em}
\mathop{\leqslant}_{\LC\RC}^{\SG_l}\nolimits\hspace{0.05em} \s'$
and $\t \hspace{0.1em}
\mathop{\leqslant}_{\LC\RC}^{\SG_{l,n-l}}\nolimits\hspace{0.1em}\t'$. 
Moreover, $\s_l \s' \hspace{0.1em}\mathop{\leqslant}_{\LC\RC}^{\SG_l}\nolimits
\hspace{0.1em} \s_l \t'$. 
Therefore, since Lusztig's conjecture $P_4$ holds in the symmetric groups, 
we obtain (a). 

If moreover $\alpb(z)=\alpb(z')$, then $\ab_{\SG_l}(\s)=\ab_{\SG_l}(\s')$ 
so $\s \sim_{\LC\RC} \s'$ by property $P_{11}$ for the symmetric group. 
Similarly, $\t \sim_{\LC\RC} \t'$ so $\s_z \sim_{\LC\RC} \s_{z'}$. 
So, by Corollary \ref{cellules bilateres}, $z \sim_{\LC\RC} z'$.\fin

\bigskip

The next proposition relates the functions $\alpb$ and $\D$. 

\medskip

\begin{prop}\label{presque P1}
Let $z \in W$. Then $\alpb(z) \le \D(z)$. Moreover, $\alpb(z)=\D(z)$ 
if and only if $z^2=1$.
\end{prop}

\proof Let us start with two results concerning the degree of $\t(\G_z)$ for 
$z \in W_n$~:

\begin{lem}\label{tau gamma}
Let $z \in W_n$. 
Then 
$$\t(\G_z)=\begin{cases}
           0 & \text{if } a_z \not= b_z, \\
	   V^{-\ell_t(z)}\t(T_{\s_{\ell_t(z)}}^{-1} C_{\s_z}) & \text{if }a_z=b_z.
	   \end{cases}$$
\end{lem}

\smallskip

\noindent{\sc Proof of Lemma \ref{tau gamma} - } 
Write $l=\ell_t(z)$. Then, $\t(\G_z)=\t(\Pi_0^T(\G_z))$. So, 
by Proposition \ref{ideaux} (a) and Corollary \ref{tau cal}, we have 
$\t(\G_z)=V^{-l}\t(T_{a_z} T_{\s_l}^{-1} C_{\s_z} T_{b_z^{-1}})$. 
Therefore, $V^l\t(\G_z)$ is equal to the coefficient of 
$T_{b_z}$ in $T_{a_z} T_{\s_l}^{-1} C_{\s_z}$. Write 
$T_{\s_l}^{-1} C_{\s_z}=\sum_{x \in \SG_{l,n-l}} \b_x T_x$. Then 
$T_{a_z} T_{\s_l}^{-1} C_{\s_z}=\sum_{x \in \SG_{l,n-l}} \b_x T_{a_z x}$. 
Thus, if $a_z \not= b_z$, then $b_z \not\in a_z \SG_{l,n-l}$ so 
$\t(\G_z)=0$. If $a_z = b_z$, then 
$\t(\G_z)=V^{-l} \b_1=V^{-l} \t(T_{\s_l}^{-1} C_{\s_z})$.\fin

\medskip

\begin{coro}\label{coro gamma}
Let $z \in W_n$. Then~:
\begin{itemize}
\itemth{a} $\deg \t(\G_z) \le -\alpb(z)$.

\itemth{b} $\deg \t(\G_z)=-\alpb(z)$ if and only if $z$ is an involution.
\end{itemize}
\end{coro}

\noindent{\sc Proof of Corollary \ref{coro gamma} - } 
This follows from Lemma \ref{tau gamma} and Corollary \ref{tw0 cw} 
(recall that Lusztig's conjectures $(P_i)_{1 \le i \le 15}$ hold 
in the symmetric group).\fin

\medskip

Let us now come back to the computation of $\D(z)$. By Corollary 
\ref{geck formule}, we have
$$\t(C_z)=\t(\G_z) + \sum_{y \inferieur z} \pi_{y,z}^* \t(\G_y).$$
But, if $y \inferieur z$, then $\alpb(z) \le \alpb(y)$ (see Proposition 
\ref{decroissante} (a)). Therefore, by Corollary \ref{coro gamma} (a), 
we have $\deg \pi_{y,z}^* \t(\G_y) < - \alpb(z)$. So 
$\deg \t(C_z) \le -\alpb(z)$ and $\deg \t(C_z)=-\alpb(z)$ 
if and only if $\deg \t(\G_z)=-\alpb(z)$ that is, if and only 
if $z$ is an involution (see Corollary \ref{coro gamma} (b)).\fin

\medskip

\section{Specialization\label{specialization}}

\medskip

We fix now a totally ordered abelian group $\G^\circ$ and a weight function 
$L^\circ : W_n \to \G^\circ$ 
such that $L^\circ(s) > 0$ for every $s \in S_n$. Let $A^\circ=\ZM[\G^\circ]$ 
be denoted exponentially and let $\HC_n^\circ=\HC(W_n,S_n,L^\circ)$. 
Let $(T_w^\circ)_{w \in W_n}$ denote the usual $A^\circ$-basis of $\HC_n^\circ$ 
and let $(C_w^\circ)_{w \in W_n}$ denote the Kazhdan-Lusztig basis of 
$\HC_n^\circ$. 

Let $b=L^\circ(t)$ and $a=L^\circ(s_1)=\dots=L^\circ(s_{n-1})$. 
Let $\th_\G : \G \to \G^\circ$, $(r,s) \mapsto ar+bs$. It is a morphism 
of groups which induces a morphism of $\ZM$-algebras 
$\th_A : A \to A^\circ$ such 
that $\th_A(V)=v^b$ and $\th_A(v)=v^a$. If $\HC_n^\circ$ is viewed as an $A$-algebra 
through $\th_A$, then there is a unique morphism of 
$A$-algebras $\th_\HC : \HC_n \to \HC_n^\circ$ such that 
$\th_\HC(T_w)=T_w^\circ$ for every $w \in W_n$. 
The main result of this section is the following:

\medskip

\begin{prop}\label{specialisation}
If $b > (n-1) a$, then $\th_\HC(C_w)=C_w^\circ$ for every $w \in W_n$.
\end{prop}

\proof Assume that $b > (n-1) a$. 
Since $\overline{\th_\HC(C_w)}=\th_\HC(C_w)$, it is sufficient to show 
that $\th_\HC(C_w) \in T_w^\circ + (\oplus_{y < w} A_{< 0}^\circ T_y^\circ)$. 
Since $\th_A(\pi_{y,w}^*) \in A_{< 0}^\circ$ for every $y < w$, 
it is sufficient to show that 
$\th_\HC(\G_w) \in T_w^\circ + (\oplus_{y < w} A_{< 0}^\circ T_y^\circ)$. 
For simplification, we set $l=\ell_t(w)$, $a = a_w$, $b=b_w$ and $\s=\s_w$. 
We set $\G_w' = T_a C_{a_l} T_\s T_{b^{-1}}$. Then 
$\G_w = \sum_{\t \le \s} p_{\t,\s}^* \G_{a a_l \t b^{-1}}'$, with 
$p_{\t,\s}^* \in v^{-1} \ZM[v^{-1}]$ if $\t < \s$ and $p_{\s,\s}^*=1$. 
So it is sufficient to show that 
$\th_\HC(\G_w') \in T_w^\circ + (\oplus_{y < w} A_{< 0}^\circ T_y^\circ)$.
By Proposition \ref{cal pl}, we have 
\eqna
\G_w' &=& T_a P_l T_{\s_l}^{-1} T_\s T_{b^{-1}}\\
&=& \DS{\sum_{k=0}^l}V^{k-l}\Big(\DS{\sum_{1 \le i_1 < i_2 < \dots < i_k \le l}} 
T_a T_{t_{i_1}t_{i_2} \dots t_{i_k}} T_{\s_l}^{-1} T_{\s b^{-1}}\Bigr)\\
&=& \DS{\sum_{k=0}^l}V^{k-l}\Big(\DS{\sum_{1 \le i_1 < i_2 < \dots < i_k \le l}} 
T_{a\a(i_1,\dots,i_k)} T_{a_k} T_{\b(i_1,\dots,i_k)}T_{\s_l}^{-1} 
T_{\s b^{-1}} \Bigr),\\
\endeqna
where $t_{i_1} \dots t_{i_k}$ is equal to $\a(i_1,\dots,i_k) a_k \b(i_1,\dots,i_k)$ with 
$\a(i_1,\dots,i_k) \in Y_{k,n-k} \cap \SG_l$ and $\b(i_1,\dots,i_k) \in \SG_l$. 
Note that $\a(i_1,\dots,i_k) a_k = r_{i_1} \dots r_{i_k}$ (recall that $r_i$ 
is defined as in \cite[\SEC 4.1]{lacriced}) so that 
$\ell(\b(i_1,\dots,i_k)) = (i_1-1) + \dots + (i_k -1)$. 
Now, let $\g(i_1,\dots,i_k)=\s_l\b(i_1,\dots,i_k)^{-1}$. 
Then 
$$\G_w'=T_w + 
\DS{\sum_{k=0}^{l-1}}V^{k-l}\Big(\DS{\sum_{1 \le i_1 < i_2 < \dots < i_k \le l}} 
T_{a \a(i_1,\dots,i_k)} T_{a_k} T_{\g(i_1,\dots,i_k)}^{-1} T_{\s b^{-1}} \Bigr).$$
If $0 \le k \le l-1 \le n-1$, we define 
$$Y_{k,l-k,n-l}=\{\s \in \SG_n~|~\forall~i \in \{1,2,\dots,n-1\} \setminus\{k,l\},~
\s s_i > \s\}.$$
Then $Y_{k,l-k,n-l}=Y_{l,n-l} (Y_{k,n-k} \cap \SG_l)$. Therefore, 
$a \a(i_1,\dots,i_k) \in Y_{k,l-k,n-l}$. But, we have 
$Y_{k,l-k,n-l} = Y_{k,n-k} (Y_{l,n-l} \cap \SG_{k,n-k})$. So we can write 
$a \a(i_1,\dots,i_k)=\a_{i_1,\dots,i_k} \a'(i_1,\dots,i_k)$ with 
$\a_{i_1,\dots,i_k} \in Y_{k,n-k}$ and 
$\a'(i_1,\dots,i_k) \in Y_{l,n-l} \cap \SG_{k,n-k}$. Then 
$\ell(\a'(i_1,\dots,i_k)) \le (l-k)(n-l)$ (indeed, 
$Y_{l,n-l} \cap \SG_{k,n-k}$ may be identified with the set of minimal 
length coset representatives of $\SG_{n-k}/\SG_{l-k,n-l}$). 
Note also that $a_k$ and $\a'(i_1,\dots,i_k)$ commute. So 
$$\G_w'=T_w + \DS{\sum_{k=0}^{l-1}}V^{k-l}
\Bigl(\DS{\sum_{1 \le i_1 < i_2 < \dots < i_k \le l}} 
T_{\a_{i_1,\dots,i_k}a_k} T_{\a'(i_1,\dots,i_k)} 
T_{\g(i_1,\dots,i_k)}^{-1} T_{\s b^{-1}} \Bigr).$$
If we write $T_u T_v^{-1} T_{\s b^{-1}} = \sum_{\t \in \SG_n} \eta_{u,v,\t} T_\t$ 
with $\eta_{u,v,\t} \in \ZM[v,v^{-1}]$, then, by \cite[Lemma 10.4 (c)]{lusztig}, 
we have $\deg \eta_{u,v,\t} \le \ell(u)+\ell(v)$. Moreover, 
$$\G_w' = T_w + \DS{\sum_{k=0}^{l-1}}V^{k-l}
\Big(\DS{\sum_{1 \le i_1 < i_2 < \dots < i_k \le l}} 
\bigl(\sum_{\t \in \SG_n} \eta_{\a'(i_1,\dots,i_k),\g(i_1,\dots,i_k),\t} 
T_{\a_{i_1,\dots,i_k} a_k \t}\bigr)\Bigr).$$
So it is sufficient to show that, for every $k \in \{0,1,\dots,l-1\}$ 
and every sequence $1 \le i_1 < \dots < i_k \le l$, we have 
$$(k-l)b + \Bigl(\ell( \a'(i_1,\dots,i_k)) + \ell(\g(i_1,\dots,i_k))\Bigr)a < 0.
\leqno{(*)}$$
But, $\ell(\a'(i_1,\dots,i_k)) \le (l-k)(n-l)$ and 
\eqna
\ell(\g(i_1,\dots,i_k)) &=& \ell(\s_l)-\ell(\b(i_1,\dots,i_k)) \\
&=& \DS{ \frac{l(l-1)}{2} } - (i_1-1) - \dots - (i_k-1) \\
&\le & \DS{\frac{l(l-1)}{2}-\frac{k(k-1)}{2}} \\
&=& \DS{\frac{1}{2}} (l-k)(l+k-1).
\endeqna
So, in order to prove $(*)$, it is sufficient to prove that 
$$2(k-l)b +  (l-k)(2(n-l) + (l+k-1))a < 0.\leqno{(**)}$$
But, 
$$2(k-l)b + a (l-k)(2(n-l) + (l+k-1))a = 
2(k-l)(b-(n-1)a) + (l-k)(k+1-l)a.$$
Since $k-l < 0$, $b-(n-1)a > 0$ and $k+1-l \le 0$, we get $(**)$.\fin

\bigskip

If $x$ and $y$ are two elements of $W_n$, we write
$$C_x^\circ C_y^\circ = \sum_{z \in W_n} h_{x,y,z}^\circ C_z^\circ,$$
where $h_{x,y,z}^\circ \in A^\circ$. We denote by $\lelo$, $\lero$, 
$\lelro$ the preorders $\lel$, $\ler$ and $\lelr$ defined in 
$\HC_n^\circ$. Similarly, we define $\sim_\LC^\circ$, $\sim_\RC^\circ$ 
and $\sim_{\LC\RC}^\circ$.

\medskip

\begin{coro}\label{cellules specialization}
Assume that $b > (n-1) a$. Let $x$, $y$ and $z$ be elements of 
$W_n$ and let $? \in \{\LC,\RC,\LC\RC\}$. Then:
\begin{itemize}
\itemth{a} $h_{x,y,z}^\circ=\th_A(h_{x,y,z})$.

\itemth{b} If $x \lepointo y$, then $x \lepoint y$.

\itemth{c} $x \sim_?^\circ y$ if and only if $x \sim_? y$.
\end{itemize}
\end{coro}

\proof (a) follows from Proposition \ref{specialisation}. (b) follows 
from (a). (c) follows from (b) and from the counting argument 
in the proof of \cite[Theorem 7.7]{lacriced}.\fin

\bigskip

Let $\t^\circ : \HC_n^\circ \to A^\circ$ denote the canonical 
symmetrizing form. If $z \in W_n$, we set
$$\ab^\circ(z) = \max_{x,y \in W_n} \deg h_{x,y,z}^\circ,$$
$$\D^\circ(z) = - \deg \t^\circ(C_z^\circ)$$
$$\alpb^\circ(z) = \th_\G(\alpb(z)).\leqno{\text{and}}$$
By Corollary \ref{cellules specialization} (b) and by 
the same argument as in the proof of Proposition \ref{decroissante}, we 
have, for every $z$, $z' \in W_n$ such that $\ell_t(z)=\ell_t(z')$ and 
$z \lelro z'$, 
\equat\label{decroissante zero}
\alpb^\circ(z') \le \alpb^\circ(z').
\endequat

\medskip

\noindent{Remark - } 
Using the result of this section, Geck and Iancu \cite{lacrimeinolf} 
proved that $\ab^\circ=\alpb^\circ$ whenever $b > (n-1) a$.\finl

\medskip

\begin{prop}\label{a special}
Assume that $b > (n-1) a$. Let $z \in W_n$. Then:
\begin{itemize}
\itemth{a} $\D^\circ(z) = \th_\G(\D(z)) \ge \alpb^\circ (z)$.

\itemth{b} $\D^\circ(z)=\alpb^\circ(z)$ if and only if $z^2=1$.
\end{itemize}
\end{prop}

\medskip

\proof First, note that $\t^\circ \circ \th_\HC = \th_\HC \circ \t$. Moreover, 
by Proposition \ref{specialisation}, we have $\th_\HC(C_z)=C_z^\circ$. 
Since $V^{\ell_t(z)} \t(C_z) \in \ZM[v,v^{-1}]$, we get that 
$\D^\circ(z) = \th_\G(\D(z))$. The other assertions 
follow easily.\fin

\medskip

We conclude this section by showing that the bound given 
by Proposition \ref{specialisation} is optimal. 

\medskip

\begin{prop}\label{optimal}
If $b \le (n-1) a$, there exists $w \in W_n$ such that 
$\th_\HC(C_w) \not= C_w^\circ$.
\end{prop}

\proof Assume that $b \le (n-1) a$. To prove the proposition, 
it is sufficient to show that there exists $w \in W_n$ such that 
$\th_\HC(C_w) \not\in T_w^\circ + \oplus_{y < w} A_{<0}^\circ T_y^\circ$. 
Using Corollary \ref{geck formule}, we see that it is sufficient 
to show that there exists $w \in W_n$ such that 
$\th_\HC(\G_w) \not\in T_w^\circ + \oplus_{y < w} A_{<0}^\circ T_y^\circ$. 
This follows from the next lemma:

\begin{quotation}
\begin{lem}\label{optimal gamma}
Let $w=s_{n-1}\dots s_2 s_1 t \s_n$. Then 
$\th_\HC(\G_w) \not\in T_w^\circ + \oplus_{y < w} A_{<0}^\circ T_y^\circ$.
\end{lem}

\smallskip

\proof We have, by Proposition \ref{cal csigma}, 
$$\G_w= T_{s_{n-1}\dots s_2 s_1 t} C_{\s_n} + 
V^{-1} T_{s_{n-1}\dots s_2 s_1} C_{\s_n}.$$
But, $T_{s_{n-1}\dots s_2 s_1} C_{\s_n}=v^{n-1} C_{\s_n}$ (see 
\cite[Theorem 6.6 (b)]{lusztig}). Therefore, 
since $\th_\HC(C_\s)=C_\s^\circ$ for every $\s \in \SG_n$, we have 
$$\th_\HC(\G_w) = \Bigl(\sum_{\t \in \SG_n} v^{(\ell(\t)-\ell(\s_n))a} 
T_{s_{n-1}\dots s_2 s_1 t \t}^\circ\Bigr) + 
v^{-b+(n-1)a} C_{\s_n}^\circ.$$
(Recall that $C_{\s_n}=\sum_{\t \in \SG_n} v^{(\ell(\t)-\ell(\s_n)} T_\t$ 
by \cite[Corollary 12.2]{lusztig}.) 
So the coefficient of $\th_\HC(\G_w)$ on $T_{\s_n}^\circ$ is 
equal to $v^{-b+(n-1)a}$, which does not belong to $A_{<0}^\circ$.\fin
\end{quotation}

\end{document}